\documentclass{article}

\usepackage{amsmath}
\usepackage[francais]{babel}
\usepackage{amssymb}
\usepackage[applemac]{inputenc}

\def\n{\noindent}
\def\Q{\mathbb Q}
\def\F{\mathbb F}

\def\Z{\mathbb Z}

\def\R{\mathbb R}
\def\C{\mathbb C}
\def\U{\mathbb U}

\def\deg{\mathop{\rm deg}}

\def\dim{\mathop{\rm dim}}

\def\a{\alpha^{\sf{v}}}
\def\rau{\rho^{\sf{v}}}
\def\er{R^{\sf{v}}}

\def\det{\mathop{\rm det}}
\def\Ad{\mathop{\rm Ad}}
\def\Tr{\mathop{\rm tr}}
\def\tr{\mathop{\rm tr}}

\def\<{\mathop{\langle}}
\def\>{\mathop{\rangle}}
\def \l{\mathop \ell}

\def\GL{\mathop{\rm GL}}
\def\SL{\mathop{\rm SL}}

\def\G{\mathop{\rm G}}
\def\Gm{\mathop{\rm G_m}}

\def\PGL{\mathop{\rm PGL}}
\def\PSL{\mathop{\rm PSL}}
\def\U{\mathop{\rm U}}
\def\SU{\mathop{\rm SU}}

\def\Lie{\mathop{\rm Lie}}
\def\Aut{\mathop{\rm Aut}}
\def\Hom{\mathop{\rm Hom}}
\def\Im{\mathop{\rm Im}}

\def\l{\ell}

\begin{document}

 \hspace{3.8cm} {\bf Zéros de caractères}

\vspace{0.3cm}

 \centerline{Jean-Pierre Serre} 
 
\bigskip

\n {\bf Introduction}

\medskip 
  Dans les $\S\S$ 1-3, nous étendons aux groupes algébriques en caractéristique zéro, et aux groupes de Lie réels compacts, le théorème de Burnside sur les groupes finis disant que {\it tout caractère irréductible de degré $>1$ s'annule au moins une fois}. En fait, lorsque le groupe n'est pas fini,
  on demande que le caractère s'annule sur un {\it élément d'ordre fini} du groupe, cf. th. 1.3 et th. 3.2.

Ce résultat avait été énoncé dans un exposé à Oberwolfach ([Se 04]), ainsi que dans un colloque en mémoire de Jacques Tits (Collège de France, décembre 2023); la démonstration, qui date de 1996, n'avait pas été publiée.

Nous donnons aussi (th. 2.16) un résultat récent de Deligne disant que,
si le groupe est  connexe, on peut renforcer la condition {\it élément d'ordre fini} en {\it élément d'ordre une puissance d'un nombre premier}. L'extension
éventuelle de ce résultat aux groupes non connexes est discutée au  $\S$ 4.8, en relation avec la classification des groupes finis simples, cf. Question 4.10.

Les $\S\S$ 4-5  portent sur les zéros des $S$-caractères (au sens de E.M. Zhmud) des groupes de Lie réels compacts.

Le $\S$ 6 est relatif au nombre (fini ou infini) de classes de conjugaison d'éléments d'ordre fini en lesquels s'annule un caractère donné. Le $\S$  7 traite l'exemple du caractère de la représentation adjointe d'un groupe de type $\sf G$$_2$: ce caractère s'annule sur des éléments d'ordre 7, 8, 15, 42.

\smallskip

J'ai plaisir à remercier C. Smyth, P. Deligne et G. Malle pour  leurs contributions au $\S$ 7, et aux $\S\S$ 2, 4.

\bigskip

\n  {\bf $\S$ 1. Le cas des groupes algébriques - préliminaires.}

\medskip

\n {\bf 1.1. Rappels sur les représentations linéaires et leurs caractères.}

Soit $k$ un corps algébriquement clos de caractéristique $0$, et soit $G$ un $k$-groupe algébrique linéaire. Le groupe $G(k)$ de ses $k$-points sera
écrit ``$G$'' dans la suite; la composante connexe de l'élément neutre de $G$
sera notée $G^0$.

 Une représentation linéaire  de $G$ (ou un ``$G$-module'') est un $k$-espace vectoriel $V$ de dimension finie muni d'un morphisme de groupes algébriques $\rho : G \to \GL_V$, où $\GL_V$ est le groupe des automorphismes de $V$, considéré comme groupe algébrique. 

Le {\it caractère} du $G$-module $V$ est la fonction $g \mapsto \tr(\rho(g))$. C'est un élément de l'algèbre affine de $G$. Elle ne dépend que des quotients de Jordan-Hölder de $V$. Deux représentations semi-simples qui ont le même caractère sont isomorphes. Un caractère est dit irréductible si un (donc tout) $G$-module  semi-simple correspondant est simple. 

Le {\it degré} d'un caractère $\chi$ est l'entier $\chi(1)$, autrement dit la dimension d'un $G$-module correspondant.

\medskip

\n {\bf 1.2. Enoncé du théorème.}

\medskip

\n {\bf Théorème 1.3.} {\it Soit $G$ un $k$-groupe algébrique linéaire et soit $\chi$ un caractère irréductible de $G$ de degré $>1$. Il existe un élément $g$ de  $G$, d'ordre fini, tel que $\chi(g)=0$.}

\smallskip

La démonstration sera donnée au $\S$ 2.14. 

\smallskip

\n {\it Remarque.} Noter la condition ``$g$ d'ordre fini''. Si on l'omet, l'énoncé devient plus facile : lorsque $G$ est connexe, il est alors valable pour tout caractère (irréductible ou non) qui n'est pas un multiple d'un caractère de degré 1.

\n{\small [Esquisse de démonstration : on peut supposer que $G$ est réductif; le cas où c'est un tore de dimension 1 résulte de ce que tout polynôme de Laurent qui n'est pas un monôme a un zéro; le cas où c'est un tore de dimension quelconque se ramène au cas de dimension 1 par restriction à des sous-tores convenables; le cas général se déduit du cas précédent, appliqué à un tore maximal de $G$.]}
\medskip

\n {\bf 1.4. Relèvements d'éléments et de sous-groupes finis.}

\medskip
\n {\bf Proposition 1.5.} {\it Soit $\phi : H_1\to H_2$ un homomorphisme surjectif de $k$-groupes algébriques linéaires, et soit $F_2$ un sous-groupe fini de $H_2$.
Il existe un sous-groupe fini $F_1$ de $H_1$ tel que $\phi(F_1)=F_2.$}

\smallskip \n  {\it Démonstration.}  Soit $G_1= \phi^{-1}(F_2)$. D'après 
le lemme 5.11 de [BSe 64]  il existe un sous-groupe fini $F_1$ de $G_1$
tel que $G_1 = G^0_1\!\cdot\!F_1$, où $G^0_1$ est la composante neutre de $G_1$, cf. 1.1. Le groupe $\phi(G^0_1)$ est un sous-groupe connexe de $F_2$; il est donc trivial. Comme $\phi(G_1)=F_2$, cela entraîne $\phi(F_1)=F_2$.

\smallskip
\n {\bf Corollaire 1.6}. {\it Tout élément d'ordre fini de $H_2$ est l'image par $\phi$ d'un élément d'ordre fini de $H_1$.}

\smallskip
 L'une des conséquences du cor. 1.6 est que, si le théorème 1.3 est vrai pour un
 groupe, il l'est aussi pour tout quotient de ce groupe.

\medskip

\n {\bf 1.7. Réductions.}

\smallskip
\n {\bf Lemme 1.8}. {\it Si le théorème 1.3 est vrai pour $G^0$, il est vrai pour $G$.}

\smallskip
\n {\it Démonstration}. D'après la prop. 1.5, il existe un sous-groupe fini $F$ de $G$ tel que $G = G^0 \!  \cdot \! F$. Soit $F$ un tel sous-groupe. Soit $V$ une représentation semi-simple de $G$ de caractère $\chi$.
Comme $G^0$ est normal dans $G$, le $G^0$-module $V$ est semi-simple
et se décompose en composantes isotypiques  $$V \ = \ \oplus_{\lambda} \ V_\lambda,$$ où les $\lambda$ sont les caractères irréductibles de $G^0$;
chaque $V_\lambda$ est somme directe de représentations irréductibles de
$G^0$ de caractère $\lambda$. Soit $L$ l'ensemble des $\lambda$ tels que
$V_\lambda \neq 0$. Le groupe $G/G^0$ opère sur $L$; cette action est
transitive, puisque $V$ est irréductible; il en est donc de même de l'action de $F$. 

Si $|L| > 1$, un théorème de Jordan
([Jo 72], [CC 92]) montre qu'il existe $g\in F$ qui opère sans point fixe sur $L$. Cela entraîne que la trace de $g$ est 0; comme $g$ appartient à $F$, il est d'ordre fini.

Lorsque $L$ a un seul élément $\lambda$, on a
$V= V_\lambda$.  Si $\deg(\lambda) > 1$, l'hypothèse faite sur $G^0$ entraîne qu'il existe $g\in G^0$, d'ordre fini,  tel
que $\lambda(g)=0$, d'où $\chi(g)=0$ puisque $\chi(g)$ est un multiple de $\lambda(g)$. Si $\deg(\lambda)=1$, $\lambda$ est un homomorphisme de $G^0$
dans le groupe multiplicatif $\Gm = \GL_1$ et $G^0$ opère sur $V$ par des homothéties; tout sous-espace de $V$ stable par $F$ est donc stable par $G$. Il en résulte
que $V$ est $F$-irréductible. En appliquant à $F$ le théorème de Burnside cité dans l'introduction, on voit qu'il existe $g\in F$ tel que $\chi(g)=0$; puisque $g$ appartient à $F$, $g$ est d'ordre fini.

\smallskip
\n {\bf Lemme 1.9}. {\it Si le théorème 1.3 est vrai pour des groupes $G_i$, il est vrai pour le groupe $G= \prod G_i$.}  

\smallskip
\n {\it Démonstration}. Le caractère $\chi$ de $G$ est produit de caractères 
irréductibles $\chi_i$ des $G_i$; comme $\deg(\chi) > 1$, il existe un indice $j$ tel que $\deg(\chi_j)>1$, et il existe donc $g_j\in G_j$, d'ordre fini, tel que $\chi_j(g_j)=0$. L'élément $g$ de $G$ dont la $j$-ième composante est $g_j$ et dont les autres composantes sont égales à 1 est d'ordre fini, et l'on a $\chi(g)=0$.

\medskip

\n {\bf Lemme 1.10}. {\it Si le théorème 1.3 est vrai lorsque $G$ est connexe, presque simple et simplement connexe, il est vrai pour tout  $G$.}

\n [Un groupe semi-simple est dit presque simple si son algèbre de Lie est une algèbre simple.]

\smallskip  

\n {\it Démonstration.} D'après le lemme 1.8, nous pouvons supposer que $G$ est connexe; nous pouvons aussi supposer que la représentation linéaire $\rho: G \to \GL_n$ correspondant à $\chi $ est fidèle (si elle ne l'est pas, on remplace $G$ par son image $\rho(G)$ - c'est licite, car tout élément d'ordre fini de $\rho(G)$ est l'image par $\rho$ d'un élément d'ordre fini de $G$, cf. cor. 1.6). 

Comme $\rho$ est irréductible, donc semi-simple, et fidèle,  le groupe $G$ est réductif connexe. Il existe alors un revêtement fini $G'$ de $G$ de la forme $T \times G_1 $, où $T$ est un tore (i.e. un produit de copies de $\Gm$), et  où $G_1$ est semi-simple et simplement connexe. Soit $\rho'$ le composé : $G' \to G \stackrel{\rho} \rightarrow\GL_n$.   Par hypothèse, le th. 1.3 est vrai pour $G_1$; il est vrai pour $T$
puisque les représentations irréductibles de $T$ sont de degré 1. D'après le lemme 1.9, il est vrai pour $G'=T \times G_1$. Il existe donc un élément 
$g'$ de $G'$ d'ordre fini tel que $\Tr \rho'(g')=0$; l'image $g$ de $g'$ dans $G$ est d'ordre fini, et l'on a $\chi(g)=0$.

\medskip
\n  { \bf $\S$ 2. Rappels sur les $\SL_2$ principaux, et fin de la démonstration du théorème 1.3}.

\medskip
  Dans ce $\S$ on suppose que $G$ possède les propriétés énumérées dans le lemme 1.10, autrement dit que c'est un groupe réductif  connexe presque simple et simplement connexe. Si l'on démontre le th. 1.3 dans ce cas, le lemme 1.10 entraîne qu'il est vrai dans le cas général.
  
  \bigskip
\n   {\bf 2.1. Poids, racines, caractères.} 

\smallskip
  Rappelons quelques notations standards ([Ja 03], II, chap. 1).
  
  \smallskip
  
 Soit $T$ un tore maximal de $G$, et soit $X = \Hom(T, \Gm)$  le groupe des caractères  de $T$ (on note $\Gm$ le tore standard de dimension 1, i.e. $\GL_1$). Soit $R \subset X$ le système de racines associé à $(G,T)$, autrement dit l'ensemble des poids $\neq 0$ de la représentation
naturelle de $T$ sur l'algèbre de Lie de $G$. Comme $G$ est presque simple,
$R$ est irréductible, au sens de [Bo 68], VI.1.2. Soit $N$ le normalisateur de $T$ dans $G$. Le quotient $W = N/T$ est le groupe de Weyl;  il opère sur $T, X$ et $R$.

 Choisissons une base $B$ de $R$, au sens de [Bo 68], VI.1.5; cela revient à choisir un sous-groupe de Borel de $G$ contenant $T$. Tout élément $\alpha$ de $R$ est combinaison $\Z$-linéaire des éléments de $B$, à coefficients, soit tous $\geqslant 0$ (auquel cas on écrit $\alpha >0$) , soit tous $\leqslant 0$. 

Notons $\er$ l'ensemble des racines $\a$ inverses des $\alpha$ (cf. [Bo 68], VI.1.1 et [Ja 03], II.1.3); c'est une partie du $\Z$-dual $Y= \Hom(\Gm, T)$ de $X$. On note
$(\alpha,\beta) \mapsto \<\alpha,\beta\>$ l'accouplement $X \times Y \to \Z$. 

\smallskip

On pose $\rho = \frac{1}{2} \sum_{\alpha > 0} \alpha$; c'est un élément de
$X$ (à cause de l'hypothèse que $G$ est simplement connexe); on a $\<\rho,\a\>=1$ pour tout $\alpha \in B$. On définit de même $\rau = \frac{1}{2} \sum_{\alpha > 0} \a$;  on a $2\rau \in Y$, mais pas toujours $\rau \in Y$.

\medskip
\n  {\bf 2.2. L'homomorphisme principal $\SL_2 \to G$.}

\medskip

Soit $f: \Gm \to T$ \ l'élément de $Y$ tel que $\<\alpha,f\> = 2$ pour tout $\alpha \in B$; on a $f = 2\rau = \sum_{\alpha > 0}\a$. 

L'élément $\varepsilon = f(-1)$ de $T$ appartient au centre de~$G$.
 On a $\varepsilon=1$ si et seulement si $f/2 =\rau$ appartient à $Y$ 
;  le noyau 
de $f$ est alors $\{1,-1\}$; c'est le cas lorsque le centre de $G$ est d'ordre impair, i.e. pour $G$ de type $\sf{A}$$_n$  ($n$ pair), ou $\sf{E}$$_6, \sf{E}$$_8, \sf{F}$$_4, \sf{G}$$_2$. 

Lorsque $\varepsilon \neq 1$,  $f$ est injectif; c'est le cas, par exemple, lorsque $G$ est
de type $\sf{C}$$_n$, ou $\sf{A}$$_n$ ($n$ impair), ou $\sf{E}$$_7$.

 ll existe
 un homomorphisme $\phi : \SL_2 \to G$ tel que le diagramme ci-dessous soit commutatif:
\vspace{1mm}

\hspace{30mm} \quad \ $\Gm \quad   \stackrel{f} \rightarrow  \quad  \ T$
\vspace{2mm}

\hspace{30mm}\quad  $\iota \downarrow  \hspace{17mm}\  \downarrow \iota_T$

\hspace{30mm}\vspace{1mm}
\quad \ \ $\SL_2\quad \stackrel{\phi} \rightarrow \hspace{5mm}  G$  

\smallskip
\n où $\iota$ désigne le plongement $t \mapsto \begin{pmatrix} t & 0 \cr 0 & t^{-1} \end{pmatrix}$ de $\Gm$ dans $\SL_2$, et $\iota_T$ désigne l'injection de $T$ dans $G$; ce résultat est dû à de Siebenthal [Si 50] et Dynkin [Dy 52].

A conjugaison près, l'homomorphisme $\phi$
ne dépend pas des choix faits; un tel homomorphisme est appelé {\it principal}; il est injectif si et seulement si l'élément $\varepsilon$ défini plus haut est $\neq 1$; lorsque $\varepsilon = 1$, $\phi$ se factorise en    $ \SL_2  \to  \PGL_2  \to  G.$

 On trouvera dans [Ko 59] de nombreuses propriétés du sous-groupe $\Im(\phi)$ de $G$.

\medskip

\n  {\bf 2.3. Calculs de caractères.}

\smallskip

 Les représentations irréductibles de $G$ sont indexées par les {\it poids dominants}, i.e. par les éléments $\lambda$ de $X$ tels que $\<\lambda,\a\>$ soit un entier  $\geqslant 0$ pour tout $\alpha >0$, cf. [St 68], chap. 12, th. 39, et [Ja 03], II, chap. 2.

  Soit $V_\lambda$ une représentation irréductible correspondant à un poids dominant $\lambda$ (autrement dit, telle que $\lambda$ soit son plus grand poids), et soit $\chi_\lambda$ le caractère correspondant de $G$. La valeur de $\chi_\lambda$ pour les éléments de $T$ est donnée par une formule, due à H. Weyl,  que nous allons rappeler.
  
   Pour écrire commodément cette formule, nous noterons exponentiellement les caractères de $T$, i.e. nous écrirons $x^\theta$ à la place de $\theta(x)$, si  $x\in T$, et $\theta\in X$.  Posons:
   
   \smallskip 
     $ A_{\theta}(x) = \sum_{w\in W} \varepsilon_wx^{w\theta},$
     où $\varepsilon_w = \det(w) \in \{1,-1\}$ est la signature de $w$.
     
     \medskip

     \n {\bf Théorème 2.4}. (H. Weyl) $\chi_\lambda(x) = A_{\lambda + \rho}(x)/A_\rho(x).$
     
     \smallskip
\n {\it Remarque.} Cette formule n'a de sens que pour les $x$ tels que $A_\rho(x) \neq 0$; mais, comme ceux-ci forment un ouvert de $T$ qui est dense pour la topologie de Zariski, elle détermine $\chi_\lambda$ sans ambiguïté. Même remarque pour la formule du th. 2.28, ainsi que pour (2.13)
ci-après.

\medskip {\it Démonstration}. Voir  [Ja 03], II, 5.11, et aussi (pour l'analogue en termes d'algèbres de Lie) [Bo 75], VIII.9.1, th.1.

   \medskip
  
Soit $y : \Gm \to T$ un élément de $Y$. Si $\theta \in X$, on a
  $y(t)^\theta = t^{\<\theta,y\>}$, d'où :
  
  \medskip
(2.5) \ $A_\theta(y(t)) = \sum_{w\in W} \varepsilon_wt^{\<w\theta,y\>}$.

\medskip
   
  \n  Le cas qui nous intéresse ici est celui où $y = f = 2\rho^{\sf{v}}$, cf. $\S$ 2.1.  
    
    \medskip
  
  \n {\bf Lemme 2.6.} {\it Pour tout $\theta \in X$, on a} $A_\theta(f(t)) = \prod_{\alpha > 0}(t^{\<\theta,\alpha^{\sf{v}}\>} - t^{-\<\theta,\alpha^{\sf{v}}\>}).$
  
  \medskip
  \n {\it Démonstration.} Soit $\Z[Y]$ la $\Z$-algèbre du groupe abélien $Y$. Cette algèbre a une base $e(y), y\in Y$, telle que
  $e(y+y') = e(y)e(y')$ pour $y, y' \in Y$. En appliquant au système de racines $R^{\sf{v}}$ une formule standard ([Bo 68], VI, 3.3, prop. 2(i)),
  on trouve : 

  \medskip
   (2.7) \  $ \sum_{w\in W} \varepsilon_we(wf) = \prod_{\alpha > 0}(e( \alpha^{\sf{v}})-  e(- \alpha^{\sf{v}}))$ \quad dans \ $\Z[Y]$.

  \medskip
  
 \n Si $\theta \in X$, soit $\phi_\theta : \Z[Y] \to \Z[t,t^{-1}]$ l'homomorphisme défini par $e(y) \mapsto t^{\<\theta,y\>}$. En appliquant $\phi_\theta$ aux deux membres de (2.7), on obtient le lemme 2.6.
 
 \medskip
  
  Nous pouvons maintenant déterminer la valeur de  $\chi_\lambda$ pour les éléments de $T$ de la forme $f(t)$ :
  
  \smallskip
    
  \n {\bf Théorème 2.8.} {\it  On a}
$$\chi_\lambda(f(t))  = \ \prod_{\alpha > 0} \frac{ t^{\<\lambda+\rho,\alpha^{\sf{v}}\>} - t^{-\<\lambda+\rho,\alpha^{\sf{v}}\>}}{t^{\<\rho,
\alpha^{\sf{v}}\>} - t^{-\<\rho,\alpha^{\sf{v}}\>}} \
\ = \ t^{- \<\lambda,2\rho^{\sf{v}}\>} \prod_{\alpha > 0} \frac{t^{2\<\lambda+\rho,\alpha^{\sf{v}}\>} - 1}{t^{2\<\rho,\alpha^{\sf{v}}\>} - 1}.$$

\medskip

[Aux notations près, le th. 2.8 se trouve dans [Ka 81], formule (10), dans le cas particulier où $t$ est une racine de l'unité (cas particulier qui est d'ailleurs suffisant pour la suite); dans le cas général, il est explicité dans [EKV 09], $\S$ 3.1, à cela près que la formule a été écrite par erreur en omettant les coefficients ``2'' du membre de droite.]

\medskip

\n  {\it Démonstration du th. 2.8.}  Il résulte du th. 2.4, combiné avec le lemme 2.6, appliqué à   $\theta = \lambda + \rho$ et à $\theta =\rho$.

\medskip

\n {\bf Corollaire 2.9}. {\it $\chi_\lambda(f(t))$ est le produit d'une puissance de $t^{-1}$ par des polynômes cyclotomiques.}

 \medskip 
 
  En effet, $\chi_\lambda(f(t))$ appartient à $\mathbb{Z}[t,t^{-1}]$ puisque c'est un caractère du groupe $\Gm$; ce polynôme de Laurent  n'a pas de pôle en dehors de $t=0$, et la formule du th. 2.8 montre que ses zéros sont des racines de l'unité.  \`{A} un facteur constant près, c'est donc
 le produit de  $ t^{- \<\lambda,2\rho^{\sf{v}}\>}$ par des polynômes cyclotomiques; on vérifie que le facteur constant est égal à 1 en observant que
 le coefficient de $ t^{- \<\lambda,2\rho^{\sf{v}}\>}$ est 1.
 
 \bigskip

\n {\it Remarque.} Lorsque l'élément $\varepsilon$ du $\S$ 2.2 est égal à 1,
$f(t)$ ne dépend que de $t^2$. Si l'on pose $u = t^2$, la formule du th.  2.8 se récrit :

  $$\chi_\lambda(f(t))  = 
 \ u^{- \<\lambda,\rho^{\sf{v}}\>} \prod_{\alpha > 0} \frac{ u^{\<\lambda+\rho,\alpha^{\sf{v}}\>} - 1}{u^{\<\rho,\alpha^{\sf{v}}\>} - 1}.$$

\bigskip

\medskip

 \n {\bf 2.10. Exemple.} 
 
 \medskip
 
 Supposons $G$ de type $\sf{G}$$_2$, ce qui entraîne $\varepsilon= 1$, et choisissons pour $\lambda$ la plus grande racine positive.  Alors
 $V_\lambda$  est la représentation adjointe de $G$, de dimension 14.  Un calcul simple montre que
 
 \medskip
   $\chi_\lambda(f(t)) = \sum_{n=-1}^{n=1} u^n +  \sum_{n=-5}^{n=5} u^n $ 
  
      \medskip
 \hspace{13mm}  $= u^5 + u^4+u^3+u^2+2u+2+2u^{-1}+u^{-2}+u^{-3}+u^{-4}+u^{-5}$
   
   \medskip
 \hspace{13mm}   $= u^{-5}\Phi_7(u)\Phi_8(u)$,

 \smallskip
\n   où $\Phi_7$ et $\Phi_8$ sont les polynômes cyclotomiques:

\medskip

  \quad \  $\Phi_7(u) = u^6+u^5+u^4+u^3+u^2+u+1, \quad \Phi_8(u)=u^4+1.$

\medskip  On obtient ainsi des éléments de $G$ d'ordre 7 et d'ordre 8 
qui sont des zéros de $\chi$. Nous verrons au $\S$ 7 qu'il y a d'autres zéros d'ordre 15 et 42.

\medskip
\n {\bf 2.11. Interprétation du théorème 2.8 en termes de produits tensoriels de représentations.}

\medskip
 \n  Introduisons d'abord quelques notations.

\smallskip
  
Pour tout $n \geqslant 1$, soit $E_n$ une représentation irréductible de $\SL_2$ de dimension $n$ (il y en a une seule, à isomorphisme près). Pour tout poids $\mu$ de $G$ tel que $\<\mu,\alpha^{\sf{v}}\> \geqslant 1$ pour tout $\alpha > 0$, posons
 $$E_\mu = \otimes_{_{\alpha > 0}} E_{\<\mu,\alpha^{\sf{v}}\>}.$$
 \n [Noter que la condition sur $\mu$ est satisfaite si et seulement si $\mu - \rho$ est un poids dominant.]

\smallskip
  Notons $\phi^*V_\lambda$ la représentation de $\SL_2$ obtenue en faisant opérer
$\SL_2$ sur $V_\lambda$ grâce à l'homomorphisme $\phi : \SL_2 \to G$. 

\medskip

\n {\bf Théorème 2.12.} {\it Les $\SL_2$-représentations  $ \phi^*V_\lambda \otimes E_{\rho} $ et $E_{\lambda + \rho}$ sont isomorphes.}
  
\medskip \n {\it Remarque}. Dans le corps des fractions de l'anneau des représentations de $\SL_2$, cet énoncé peut s'écrire : $$ \phi^*V_\lambda = E_{\lambda + \rho}\cdot (E_{\rho})^{-1}.$$ En comparant les dimensions des deux membres, on retrouve la formule
de H. Weyl donnant $\dim V_\lambda$, cf. [Bo 75], chap.  8, $\S$ 9, th. 2; c'est d'ailleurs à peu près la méthode suivie par Weyl lui-même dans [We 14].

\medskip

\n {\it Démonstration du théorème 2.12.} 

Pour prouver que deux caractères d'un groupe réductif connexe sont égaux, il suffit de le tester sur un tore maximal; en effet, les tores maximaux sont conjugués, et leur réunion est  l'ensemble des éléments semi-simples, qui est dense. 

Dans le cas de
$\SL_2$, cela signifie que ces caractères
prennent la même valeur pour les éléments  $s_t= \begin{pmatrix} t & 0 \cr 0 & t^{-1} \end{pmatrix}$ du tore maximal $ \begin{pmatrix} * & 0 \cr 0 &* \end{pmatrix}$. 

Calculons ces valeurs pour les représentations qui interviennent dans le th. 2.12 :

\medskip
Si $n\geqslant  1$, la valeur en $s_t$ du caractère de $E_n$ est 
$(t^n-t^{-n})/(t-t^{-1})$.

\medskip

 La valeur en $s_t$ du caractère de  $E(\mu)$ est $ \prod_{\alpha > 0}({t^{\<\mu,\alpha^{\sf{v}}\>} - t^{-\<\mu,\alpha^{\sf{v}}\>}})/(t-t^{-1}).$
 
 \smallskip
\n En  appliquant cette dernière formule à $\mu =\lambda + \rho$, puis à  $\mu=\rho$, et en divisant, les puissances de $t-t^{-1}$ disparaissent et les autres facteurs donnent $\chi_\lambda(f(t))$ d'après le th. 2.8. 

\smallskip

\n {\it Exemple.} Si $k = \C$, et $t = e^{i\theta}$, la valeur en $s_t$ du  caractère de $E_n$ est $\sin(n\theta)/\sin(\theta)$ et l'on a donc:

\smallskip
$$(2.13) \quad \chi_\lambda(f(e^{i\theta})) = \prod_{\alpha > 0} \frac{\sin(\<\lambda+\rho,\alpha^{\sf{v}}\>\theta)}{\sin(\<\rho,
\alpha^{\sf{v}}\> \theta)}, \quad {\rm cf. \ [Ka \ 81],  \ formule \ (10)}.$$

\smallskip

       \medskip
\n  {\bf  2.14. Fin de la démonstration du théorème 1.3.} 

\smallskip
  Conservons les notations du th.1.3, ainsi que celles du $\S$ 2.3, et notons $\lambda$ le poids dominant tel que $\chi = \chi_\lambda$. Pour  achever la démonstration du th.1.3, il suffit de montrer qu'il existe
 une racine de l'unité $z$  tel que $\chi_\lambda(f(z))=0.$
  On n'a que l'embarras du choix puisque, d'après le cor. 2.9, tous les zéros
  du polynôme de Laurent $\chi_\lambda(f(t))$ sont des racines de l'unité.
   
  On peut, par exemple, prendre $z$ d'ordre $m = \<2\lambda+2\rho,\beta^{\sf{v}}\>$, où $\beta^{\sf{v}}$ est la plus grande racine de $\er$; 
 le numérateur du membre de droite du th. 2.4 est nul, alors que le 
  dénominateur ne l'est pas, car $\<2\rho,\alpha^{\sf{v}}\> \leqslant  
 \<2\rho,\beta^{\sf{v}}\> < m$ pour tout $ \alpha> 0$. 
 
\bigskip

\n {\bf 2.15. Un théorème de Deligne.}

\smallskip
  Tout récemment, Deligne a montré que le th. 1.3 peut être grandement
  amélioré lorsque $G$ est connexe : ``ordre fini'' peut être remplacé
  par ``ordre une puissance d'un nombre premier". Autrement dit :
  
  \smallskip  
  \n {\bf Théorème 2.16.} (Deligne) {\it Soit  $\chi$ un caractère irréductible dee degré $>1$ d'un groupe algébrique connexe. Il existe un  élément $g$ du groupe, d'ordre une puissance d'un nombre premier, tel que $\chi(g)= 0$.}
  
  (Dans une version antérieure de cet article, le th. 2.16 était énoncé comme ``Problème''. )
  
\smallskip
  \n {\it Démonstration} (d'après Deligne \footnote{Voir son exposé {\it Zeros of characters}, I.A.S, Princeton, 4 novembre 2024.} - reproduite avec sa permission).
  
\smallskip

    Pour simplifier, nous dirons {\it d'ordre primaire} à la place de {\it d'ordre une puissance d'un nombre premier}.
    
\smallskip
    
    Le point essentiel est la proposition suivante:
  
\smallskip  
    \n {\bf Proposition 2.17}. {\it Soit $F(t)$ une fonction rationnelle de $t$, 
 de la forme
 $$F(t) =  t^a \prod_i(t^{n'_i} - 1)/(t^{n_i} - 1),$$
   où les $n'_i,n_i$ sont des entiers $\geqslant 1$. Supposons que 
  $\prod {n'_i} > \prod {n_i}.$ Il existe alors une racine de l'unité 
  d'ordre primaire qui est un zéro de $F$.}

\smallskip
  \n {\it Démonstration.} Pour tout nombre premier $\ell$, notons $n'_{i,\l}$ la $\l$-composante de $n'_i$ et définissons de même $n_{i,\l}$. On a
  
 $$ \prod_\ell \frac{\prod_i{n'_{i,\ell}}}{\prod_i{n_{i,\ell}}} \ = \ \prod_i \frac{n'_i}{n_i} \ > \ 1.$$
  
  Il existe donc un nombre premier $\ell$ tel que $\prod_i n'_{i,\l} > \prod_i,n_{i,\l}$.
  Choisissons un tel $\ell$. On a 
  
  $$  \sum_i v_\ell(n'_i) \  > \sum_i v_\ell(n_i),$$
  où  $v_\ell$ est la valuation $\ell$-adique de $\Q$.
  
  Pour tout entier $m > 0$,  soit $\Phi_{\ell^m}$ le polynôme cyclotomique d'indice $\ell^m$. C'est un polynôme irréductible sur $\Q$. Il définit donc une valuation discrète de $\Q(t)$, que nous noterons $w_m$. Si $n$ est un entier $> 0 $,  $w_m (t^n-1)$ est égal à 1 si $n$ est divisible par $\ell^m$, et à 0 sinon. On a donc
    $\sum_m w_m(t^n-1) = v_\ell(n),$
    d'où:
    
    $$  \sum_m w_m(F) =   \sum_i v_\ell(n'_i)  - \sum_i v_\ell(n_i) > 0.$$
    
    Il existe donc au moins un $m$ tel que $w_m(F) > 0$. Cela signifie que  $F$ a un zéro qui est racine primitive $\ell^m$-ième de l'unité.	
    
    \medskip
    
    \n {\it Fin de la démonstration du th. 2.16.}
    
    Notons d'abord que, si $G_1 \to G_2$  est un homomorphisme surjectif de groupes linéaires, tout élément de $H$ d'ordre primaire est l'image d'un élément de $G$ d'ordre primaire. Cela résulte du cor. 1.6, combiné à un
    argument à la Bézout. Cela permet de supposer que le groupe algébrique du th. 2.16 est réductif. L'argument du lemme 1.10 montre que l'on peut se ramener au cas où ce groupe est presque simple et simplement connexe. Les hypothèses du présent $\S$ sont satisfaites.
  
     On peut alors appliquer les constructions 
    des précédentes sections, et l'on obtient un homomorphisme
    $f : \G_m \to T \to G$, où $T$  est un tore maximal de $G$. Si $\lambda $ est le poids maximum du caractère $\chi$, le th. 2.8 donne une 
    formule pour $F(t) =\chi(t)$ qui est du type de la prop. 2.17, les $n'_i$
    étant les $\<\lambda+\rho,\alpha^{\sf{v}}\>$ et les $n_i$ les $\<\rho,\alpha^{\sf{v}}\>$. On a $n'_i \geqslant n_i$ pour tout $i$, avec
    inégalité stricte pour au moins un indice $i$ puisque le poids $\lambda$
    est non nul. D'où $\prod n'_i > \prod n_i$, et l'on applique la prop. 2.17;
    elle montre que $F$ a un zéro du type voulu.
    
    \medskip
    
    \n {\it Remarque.} Supposons $G$ simple. La démonstration ci-dessus montre, non seulement qu'il existe un zéro de $\chi$ qui est d'ordre primaire, mais aussi qu'un tel élément {\it peut être choisi de type principal}, i.e.  contenu dans un sous-groupe principal de type $\sf{A}$$_1$ de $G$.

\medskip
\n {\bf $\S$ 3. Le cas des groupes compacts.}

\medskip

\n  {\bf 3.1. Enoncé du théorème}.

\smallskip

  Cet énoncé est presque le même que celui du th. 1.3 :
  
  \smallskip

\n {\bf Théorème 3.2.} {\it Soit $G$ un groupe topologique compact et soit $\chi$ le caractère d'une représentation linéaire complexe continue irréductible $\rho : G \to \GL_n(\C)$ de degré $n>1$. Il existe $g\in G$ tel que $\chi(g)=0; $ si $G$ est un groupe de Lie réel, on peut choisir $g$ d'ordre fini.}

\smallskip

\n {\bf Lemme 3.3.} {\it Si le théorème 3.2 est vrai pour les groupes de Lie compacts, il est vrai pour tous les groupes compacts.}

\smallskip  \n {\it Démonstration.} Le groupe $\rho(G)$ est un sous-groupe fermé de
$\GL_n(\C)$. D'après un théorème d'Elie Cartan ([Bo 72], chap. III, $\S$ \ 8, th.2), cela entraîne que $\rho(G)$
est un groupe de  Lie réel. En appliquant le th. 3.2 à ce groupe, on obtient
un élément $\gamma \in \rho(G)$ tel que $\Tr(\gamma)=0$. Si $g\in G$
est tel que $\rho(g)=\gamma$, on a $\chi(g)=0$.

\smallskip \n {\it Remarque.} Lorsque $G$ n'est pas un groupe de Lie réel, on ne peut pas toujours choisir $g$ d'ordre fini. Voici un exemple:

  Soit $A$ le groupe additif des entiers 2-adiques, et soit $C$ un groupe cyclique d'ordre 3. Faisons opérer $A$ sur $C$ par 
  $A \to \Z/2\Z\  \simeq \Aut(C)$ et soit $G = C\!\cdot\! A$ le produit semi-direct correspondant. C'est un groupe profini, donc compact. Ses seuls éléments d'ordre fini sont ceux de $C$. Il a un quotient isomorphe au groupe symétrique $S_3$, d'où une représentation irréductible de degré 2; soit $\chi$ son caractère. Les valeurs de $\chi$ sur $C$ sont égales à 2 (élément neutre) et à $-1$ (éléments d'ordre 3). Il n'existe donc aucun $g\in G$ d'ordre fini tel que $\chi(g)=0$.

{\small Noter que le groupe $G$ ci-dessus est un groupe de Lie 2-adique de dimension 1; on ne peut donc pas supprimer le mot ``réel'' dans le  théorème 3.2.}

\medskip

\n  {\bf 3.4. Démonstration du th. 3.2.}
 
\smallskip 

D'après le lemme 3.3, on peut supposer que $G$ est un groupe de Lie réel compact. Nous allons utiliser le dictionnaire \footnote{Très bien expliqué dans l'article ``{\it Complexification $($Lie group$)$}'' de Wikipedia.} :

\smallskip

\centerline{{\it groupes de Lie réels compacts} $ \ \Longleftrightarrow $ \ {\it groupes algébriques réductifs complexes,}}

\smallskip   
\n pour nous ramener au th.1.3 :

\medskip
Soit $\underline{G}$ le $\R$-groupe algébrique dont l'algèbre affine est formée des fonctions continues $f$ sur $G$, à valeurs réelles, dont les translatées (à gauche) engendrent un $\R$-espace vectoriel de dimension finie (cf. [Ch 46], chap. VI  et [Se 93], $\S$ 5). Le groupe des $\R$-points de $\underline{G}$ est $G$. Soit $G_c$ le $\C$-groupe algébrique déduit de $\underline{G}$ par extension des scalaires; c'est un groupe linéaire réductif, appelé le {\it complexifié} de $G$. On a $\Lie(G_c) = \C \otimes_\R \Lie(G)$. Le groupe  $G $ est un sous-groupe compact maximal de $G_c$; tout sous-groupe compact de $G_c$ est conjugué d'un sous-groupe de $G$.
 
 Le caractère $\chi$ de $G$ s'étend en un caractère (algébrique) $\chi_c$ de $G_c$, qui est irréductible.
 D'après le th. 1.3, appliqué à $k = \C$ et au groupe $G_c$, il existe $g\in G_c$, d'ordre fini, tel que $\chi_c(g) = 0$. Comme le groupe engendré par $g$ est fini, donc compact, l'un de ses $G_c$-conjugués est contenu dans $G$; on obtient ainsi un zéro de $\chi$ d'ordre fini.

\bigskip

\n {\it Remarque.} La même démonstration montre que le théorème 2.16 est valable pour les groupes de Lie compact connexes:  {\it tout caractère irréductible de degré $> 1$ a un zéro qui est d'ordre primaire.}

\bigskip
\n {\bf $\S$ 4. $S$-caractères.}

\smallskip

Dans ce $\S$, les caractères sont des caractères complexes.

\medskip

  La notion de $S$-caractère a été introduite par E.M. Zhmud ([Zh 95]) 
  pour les groupes finis. Commençons par ce cas.

\medskip
\n {\bf 4.1. Les $S$-caractères d'un groupe fini.}

\smallskip
 Soit $G$ un groupe fini, et soit $f$ un caractère virtuel de $G$ (autrement dit une différence de deux caractères). On dit (Zhmud, {\it loc. cit.})
 que $f$ est un {\it $S$-caractère} s'il a les deux propriétés suivantes :

\smallskip
(a) {\it $f$ est positif, i.e. $f(g)$ est réel $\geqslant 0$ pour tout $g\in G$}. 

\smallskip
 (b) $\<f,1\>= \ 1$, {\it i.e.}  $\frac{1}{|G|} \sum_{g\in G} f(g) = 1.$
 
 \smallskip

 \n {\it Exemples de $S$-caractères} (cf. [Zh 95] et [BKZ 19], chap. XX) :
 
 \smallskip
 \n $(4.1.1)$ \ $f=1$.
 
 \n $(4.1.2)$ \ $f=1 - \varepsilon$ (resp. $f=1+ \varepsilon$), où $\varepsilon : G \to \{1,-1\}$
 est un homomorphisme surjectif.
 
 \n $(4.1.3)$ \ $f =\phi\overline{\phi}$, où $\phi$ est un caractère irréductible, et $\overline{\phi}$ est le conjugué de $\phi$.
 
 \n $(4.1.4)$ \ Le caractère associé à une action transitive de $G$ sur un ensemble non vide. 
   
  \n $(4.1.5)$ \ Le caractère induit d'un $S$-caractère d'un sous-groupe propre de $G$.
  
  \medskip
  \n {\it Remarques.}
  
  \n $(4.1.6)$ Voici un exemple de $S$-caractère qui n'est d'aucun des types ci-dessus. Soit
  $G= \PSL_2(\F_{7})$, soit $\psi$ le caractère de la représentation irréductible de dimension 8.   Les valeurs de $\psi$ sur les éléments de $G$ d'ordre $(1,2,3,4,7)$ sont $(8,0,-1,0,1)$. Le caractère $f = 1 + \psi$ est un $S$-caractère de degré 9 qui n'est évidemment pas de type (4.1.1), (4.1.2) ou (4.1.3), et pas non plus de type (4.1.4) ou (4.1.5) car son degré ne divise pas l'ordre de $G$.
  
\smallskip
\n $(4.1.7)$  Si $f= \sum_{\chi \  {\rm irr}} n_\chi \chi$ est un $S$-caractère, on a ([Zh 95], lemme 2):
   
   \medskip
\hspace{9mm}   $ -\chi(1) \leqslant n_\chi \leqslant \chi(1)$ \quad pour tout $\chi$ irréductible.
   
  \smallskip Cela montre que l'ensemble des $S$-caractères est fini, et cela permet de les énumérer à partir de la table des caractères de $G$

\medskip

 \n {\bf 4.2. Les zéros d'un $S$-caractère.}
 
 \medskip
 
 Les textes cités plus haut  donnent une série de résultats
 sur les valeurs des $S$-caractères. Nous nous bornerons à l'énoncé le plus simple :
 
 \medskip
\n  {\bf Théorème 4.2} ([Zh 95], th.1). {\it Soit $G$ un groupe fini. Si $f$ est un $S$-caractère distinct du caractère $1$, il existe $g\in G$ tel que $f(g)=0$.}
  
  \medskip
  
 \n {\it Démonstration} (inspirée de celle de Burnside ([Bu 03]). 
 
  Supposons $f(g) \neq 0$ pour tout $g\in G$. Soit $N = \prod_{g\in G} f(g)$. C'est un nombre réel $> 0$.  Nous allons voir que $N$ {\it appartient à $\Z$}.
 
  Soit $R(x,y)$ la relation d'équivalence sur $G$ définie par ``$x$ et $y$
  engendrent le même sous-groupe''. Soit $X$ une classe d'équivalence 
  de $R$ et soit $d$ l'ordre des éléments de $X$. Soit $N_X = \prod_{x\in X} f(x)$. Soit 
  $K_d$ le sous-corps de $\C$ engendré par les racines $d$-ièmes de l'unité; on a $f(x) \in K_d$ si $x\in X$. On sait (Kronecker) que le groupe de Galois de $K_d/\Q$ est canoniquement isomorphe au groupe $\Gamma_d = (\Z/d\Z)^\times$; l'isomorphisme associe à $j \in \Gamma_d$ l'automorphisme $\sigma_j$ tel que $\sigma_j (z)=z^j$ pour tout $z$ tel que $z^d=1$. L'action des $\sigma_j$ sur les $f(x), x \in X$ est :
    $\sigma_j(f(x)) = f(x^j)$, cf. par exemple [Se 71], $\S$ 12.4, th. 25.
Si $x\in X$, les  éléments de $X$ sont les puissances $x^j, j \in \Gamma_d$.
Il en résulte que $N_X$ est la norme de $f(x)$ dans l'extension $K_d/\Q$,
donc que $N_X$ appartient à $\Q$. Comme c'est un entier algébrique (car combinaison $\Z$-linéaire de racines de l'unité), c'est un élément de $\Z$. Puisque $N$ est le produit des $N_X$, on a $N \in \Z$.
   
Comme $N$ est positif et non nul, on a $N\geqslant 1$. La moyenne géométrique de la famille $\{f(g), g \in G\}$ est  $N^{1/n} \geqslant 1$, où $n =  |G|$. D'autre part la moyenne arithmétique des $f(g)$ est  $\<f,1\>= \ 1$. L'inégalité standard entre ces deux moyennes montre que $N^{1/n} \leqslant 1$, donc que les deux moyennes sont égales à $1$, ce qui entraîne $f(g)=1$ pour tout $g$, contrairement à l'hypothèse $f\neq 1$.
  
    \medskip

\n {\it Exemples}. 

Soit $\phi$ un caractère irréductible de degré $>1$. En appliquant le th. 4.2 à $f = \phi \overline{\phi}$, on retrouve le théorème de
Burnside.

  Soit $f$ le caractère associé à une action transitive de $G$ sur un  ensemble
 $X$ à au moins deux éléments:  si $g\in G$, $f(g)$ est le nombre de points de $X$  fixés par $g$. En appliquant le th. 4.2 à $f$, on voit qu'il existe un élément de $G$ qui ne fixe aucun point de $X$: on retrouve le théorème de Jordan cité dans la démonstration du lemme 1.8.

\smallskip
\n {\it Remarque.} Soit $f$ comme dans le th.4.2, soit $Z_f$ l'ensemble de
ses zéros et soit $c_f = |Z_f|/|G|$ la densité de $Z_f$. On trouvera dans
 [BKZ 19], chap. XX, diverses minorations de $c_f$. La plus simple est
celle où l'on suppose que $f$ est à valeurs dans $\Z$ :
on a alors $c_f \geqslant 1/\sup_{g\in G} f(g)$. La démonstration est  la même que celle donnée dans [CC 92] pour le cas particulier des caractères de permutation, où $\sup_{g\in G} f(g)=f(1)$.

\medskip

\n {\bf 4.3. Les $S$-caractères d'un groupe de Lie compact.}

\smallskip

  Soit $G$ un groupe de Lie réel compact. La notion de $S$-caractère se définit comme dans le cas des groupes finis : c'est un caractère virtuel $f$
  à valeurs réelles $\geqslant 0$, tel que  $\<f,1\>= \ 1$, autrement dit
  $\int_Gf(x) dx = 1$, où $dx$ est la mesure de Haar de $G$ de masse totale 1.
  
\smallskip  Cette définition suggère plusieurs problèmes. Le plus évident
  est de savoir si le th. 4.2 se généralise. Autrement dit : 
  
  \medskip
  
  \n {\bf Problème 4.4}. {\it Si $f$ est un $S$-caractère de $G$ distinct de $1$, existe-t-il $g \in G$ tel que $f(g)=0$ ? Si oui, existe-t-il un tel $g$ qui soit d'ordre fini ?}
  
  \smallskip
  
    Il est tentant d'essayer d'attaquer ce problème par la méthode de Burnside, cf. démonstration du th. 4.2. Bien sûr, si $G$
    est infini, on ne peut pas faire le produit des valeurs de $f$. Mais, si
    $f$ est réel $>0$, on peut définir l'intégrale $I = \int_G \log (f(x)) dx$.
    Comme $\int_G f(x) dx = 1$, l'analogue en théorie de l'intégration de l'inégalité de la moyenne ({\it inégalité de Jensen})
    montre que  $I$ est $ \leqslant 0$, avec égalité seulement si $f=1$. Pour conclure, il faudrait pouvoir prouver que $I$ est $\geqslant 0$; lorsque $G$ est fini, c'est un argument de nature galoisienne qui a permis de le faire; par quoi peut-on le remplacer ?

      \medskip
  Le second problème reflète la difficulté qu'il y a à construire des $S$-caractères lorsque $G$ est connexe :
  
  \smallskip
  
  \n {\bf Problème 4.5.} {\it Si $f$ est un $S$-caractère de $G$, et si $G$ est connexe, est-il vrai que $f$ est de la forme $f =\phi\overline{\phi}$, où $\phi$ est un caractère irréductible d'un revêtement connexe fini de  $G$ ?}
  
  \smallskip
  
Précisons ce que signifie la fin de l'énoncé: soit $G'$ un revêtement connexe fini de $G$, et soit $N$ le noyau de la projection $G'\to G$; le groupe $N$ est un sous-groupe fini du centre de $G'$. Si $\phi$ est un caractère de $G'$, correspondant à une représentation irréductible $\rho : G' \to \GL_n(\C)$, l'image de $N$ par $\rho$ est formée d'homothéties de module $1$; le produit tensoriel de $\rho$ et de sa conjuguée $\overline{\rho}$ est égal à 1 sur $N$, donc définit une représentation linéaire de $G$  dont le caractère est $\phi\overline{\phi}$; c'est évidemment un $S$-caractère.
  
  \smallskip  \n {\it Exemple. } Soit $G = \SU(n,\C)/N$, où $N$ un sous-groupe non trivial du centre de $G$ (un tel $N$ existe si $n >1$).
 Soit $f_0$
  le caractère de la représentation adjointe de $G$. Alors $f = 1+f_0$ est un $S$-caractère de $G$ de degré $n^2$; si on l'identifie à un caractère de
  $G'=\SU(n,\C)$, il est égal à $\phi\overline{\phi}$, où $\phi$ est le caractère standard de degré $n$ de  $G'$; cependant,  il n'existe pas de caractère irréductible $\psi$ de $G$ tel que $f = \psi\overline{\psi}$.
  
  \smallskip
  \n {\it Remarques.}
  
  1. Noter que le caractère $f$ ci-dessus  est de la forme $1 + \theta$, où $\theta$
  est un caractère irréductible. J'ignore s'il existe d'autres $S$-caractères $F$ qui soient de cette forme (en supposant $G$ connexe et presque simple, pour simplifier). Il n'y en a pas si le problème 4.5 a une réponse positive. En effet, supposons que 
 $F= \phi\overline{\phi}$ avec $\deg(\phi) > 1$, où $\phi$ est le caractère d'une représentation irréductible fidèle  $\rho: G \to \U(n,\C)$, et supposons que $F = 1 + \theta$, avec $\theta$ irréductible, $\theta \neq 1$. On alors
   $$\langle\phi^2\overline{\phi}^2,1\rangle = \langle F^2,1\rangle = 1 + 2\langle \theta,1\rangle + \langle \theta^2,1\rangle = 2,$$ ce qui entraîne que $G$ contient $\SU(n,\C)$ d'après {\it l'alternative de Larsen}, cf. [Bo 23] V.4.9, th. 3 a). Comme on a supposé $G$ semi-simple connexe, cela entraîne $G = \SU(n,\C)$.
  
  \smallskip
  
  2. Le th. 3.2 montre que, si le problème 4.5 a une réponse positive, il en est de même du problème 4.4 lorsque $G$ est connexe.
    
  \medskip
  
Il y a deux cas où le problème 4.5 a une réponse positive: celui où $G$ est commutatif et celui où $G= \SU(2,\C)$ :

\smallskip
  
  \n {\bf Proposition 4.6.} {\it Si $G$ est connexe et commutatif, le seul $S$-caractère de $G$ est le caractère $1$.}
  
  \smallskip
  \n {\bf Proposition 4.7.} {\it Tout $S$-caractère de $\SU(2,\C)$ est le carré 
  d'un caractère irréductible.}
  
  \medskip
  Les démonstrations seront données au $\S$ 5.5.
  
  \smallskip
  \n Ces cas sont trop particuliers pour que l'on puisse transformer les questions 4.4 et 4.5 en conjectures.
  
  \medskip
  \n{\bf 4.8. Autres problèmes sur les zéros des $S$-caractères d'un groupe de Lie compact.}
  
  \smallskip Conservons les notations du $\S$ précédent. Nous allons modifier le problème 4.4  en renforçant la condition ``$g$ est d'ordre fini''.
  Dans la première version de ce texte, le renforcement était le suivant :

    \smallskip \n {\bf Problème 4.9.}  {\it Soit $f$ un $S$-caractère de $G$ distinct de $1$. Existe-t-il un élément $g\in G$, d'ordre primaire (i.e. une puissance d'un nombre premier), tel que f$(g)=0$ ?}

    \smallskip
    C'était trop demander, même pour un groupe fini : G. Malle ([Ma 24]) a construit des contre-exemples, où $G$ est un groupe simple de type $J_1$ (ou $A_8$). Une étude détaillée a été faite ensuite par T. Breuer, M. Joswig et G. Malle ([BJM 24]).
    
    \medskip
    
    Une question plus raisonnable (surtout après le th. 2.16 de Deligne) est:
    
    \smallskip
    \n {\bf Question 4.10.} {\it Soit $\chi$ un caractère irréductible de degré $>1$ d'un groupe de Lie compact $G$. Existe-t-il un zéro de $\chi$
    qui soit d'ordre primaire?}
   
     \smallskip
      Essayons de réduire cette question au cas où $G$ est connexe.
      On procède comme dans la démonstration du lemme 8. On choisit 
      un sous-groupe fini $F$ de $G$ tel que $G = G^0 \! \cdot \! F,$
      et l'on décompose la représentation irréductible de caractère $\chi$ par rapport à l'action de $G^0$. Comme pour le lemme 1.8, on arrive ainsi à se ramener au cas $G=G^0$ traité par Deligne, mais à deux conditions sur $F$ :
      
        \smallskip
      (A) {\it Pour toute action transitive de $F$ sur un ensemble à au moins deux éléments, il existe un élément de $F$ d'ordre primaire, qui opère sans points fixes.}

      (B) {\it Pour tout caractère irréductible $\chi$ de $F$, de degré  $> 1$,
      il existe $g \in F$, d'ordre primaire, tel que $\chi(g)=0$.}
      
      \smallskip
      En fait, des réponses positives à (A) et (B) sont énoncées, et presque démontrées, dans la littérature: (A) dans  [FKS 81] et (B) dans [MNO 00].
      
       Pourquoi ``presque'' ? Parce que, dans les deux cas, la classification des groupes finis simples (CFSG) est utilisée, alors que la démonstration de ce énoncé fondamental n'est pas encore entièrement publiée.
         
  On peut espérer que, dans un futur pas trop lointain, cette
  difficulté sera levée, et que la question 4.10 sera transformée en théorème. En attendant, on peut se contenter des nombreux cas où l'on a prouvé que $F$  a les propriétés (A) et (B).
  
  \bigskip 
 
\n{\bf  $\S$ 5. Polynômes de Laurent positifs et applications aux $S$-caractères.}

\medskip
Le but de ce $\S$ est de démontrer les propositions 4.6  et 4.7.

\smallskip
\n {\bf 5.1. Quelques propriétés des polynômes de Laurent positifs.}

\medskip
Dans ce qui suit, $f(t)= \sum_{n \in \Z} a_nt^n$  est un polynôme de Laurent à coefficients réels, autrement dit un élément de l'anneau $\R[t,t^{-1}]$. On s'intéresse aux valeurs de $f$ sur le cercle unité $\U =\U(1,\C)$ de $\C$, et plus particulièrement  à la propriété de positivité :

\smallskip
  {\it  $f$ est positif, i.e. les valeurs de $f$ sur $\U$ sont des nombres réels $\geqslant 0$.}
  
  \smallskip
  Le fait que les valeurs de $f$ soient réelles  équivaut à la {\it propriété de symétrie}   $a_n=a_{-n}$ pour tout $n$. On a donc
   $$f(e^{i\theta}) = a_0 + \sum_{n \geqslant 1} 2a_n\cos(n\theta) \quad {\rm pour \ tout} \ \theta \in \R. $$
   On a $a_0 = \int_{\U} f(t)\mu_U(t)$, où $\mu_U$ est la mesure de Haar de U de masse totale 1. Si l'on écrit les éléments de U sous la forme $e^{i\theta}$, on a  $\mu_{\U} = \frac{1}{2\pi}$$d\theta$. Si $f$ est positif, on a $a_0 \geqslant 0$, avec égalité seulement si $f=0$.

\smallskip

{\bf Proposition 5.2.} {\it Supposons que $f$ soit positif. Soit $m$ un entier $ \geqslant 1$ tel  que $a_n=0$ pour tout 
$n>m$ divisible par $m$.}

\smallskip
(i) {\it \ On a}  $| a_m| \leqslant a_0/2.$

\smallskip

(ii) {\it \ Si $a_m = a_0/2$, alors $f$ est divisible par \ $t^{-m}+2+t^m$.}

\smallskip
(iii){\it \ Si  $a_m = -a_0/2$, alors $f$ est divisible par \ $-t^{-m}+2-t^m$.}

\smallskip
\n [Il s'agit de divisibilité dans l'anneau $\R[t,t^{-1}]$.]

\smallskip

\n {\it Démonstration.}

Notons $\U_+$ l'ensemble des $t\in \U$ tels que $t^m=1$, et $\U_-$ l'ensemble des $t$ tels que  $t^m=-1$. %[Autrement dit, $U_+ = \mu_m$ et $U_-= \mu_{2m} \sm \mu_m$.] 
  
  \smallskip
  Soient   $S_+= \sum_{t \in \U_+} f(t), \ S_- = \sum_{t \in \U_-} f(t).$
  
\smallskip

\n {\bf Lemme 5.3}. {\it On a $S_+ = m (a_0+ 2 a_m), \  S_- = m (a_0- 2 a_m)$ et $S_+ + S_- =  2m \!\ a_0$}.

\smallskip

\n {\it Démonstration du lemme.} Si $n$ est un entier, la somme $ \sum_{t \in \U_+} t^n$ est égale à $m$ si $m$ divise $n$, et à $0$ sinon. Les seuls termes {\it a priori} non nuls de la somme $S_+$ sont donc ceux relatifs à $n = -m, 0, m$, qui donnent 
respectivement $m \!\ a_m, m  \!\ a_0, m  \!\ a_m$; leur somme est $m (a_0+ 2 a_m)$. Le même argument, appliqué à $2m$ au lieu de $m$, montre que tous les termes de la somme $S_++S_-$ sont nuls, sauf $2ma_0$. Par différence, on en déduit  $S_- = m(a_0-2a_m)$.

\smallskip
\n {\it Démonstration de la prop. 5.2.}

Comme  $S_+$ et $S_-$ sont $\geqslant 0$, le lemme 5.3 montre que $a_0 + 2a_m \geqslant 0$ et $a_0 - 2a_m \geqslant 0$, autrement dit $- a_0/2 \leqslant a_m \leqslant a_0/2$, ce qui prouve (i). Si $a_m = a_0/2$, on a $S_-= 0$; comme $S_-$ est la somme des $f(t), t\in \U_-$, qui sont $\geqslant 0$,  tous ces $f(t)$ sont nuls.
L'ensemble $U_-$ est donc contenu dans l'ensemble des zéros de $f$. Comme $f$ est partout positive, ses zéros sont de multiplicité paire. Cela montre que $f$ est divisible par 
$(t^m+1)^2$, donc aussi par $(t^{-m}+1)(t^m+1)$, ce qui prouve (ii). La démonstration de (iii) est analogue à celle de (ii).

\medskip
\n {\bf Corollaire 5.4.} {\it Supposons que $f$ soit positif, non constant, et que l'on ait $a_n \in \Z$ pour  tout $n$.
 On a alors $a_0 \geqslant 2;$ si $a_0 = 2$, $f$ est égal, soit à $t^{-m}+2+t^m$, soit à $-t^{-m}+2-t^m$, avec $m > 0$.}
 
 \smallskip {\it Démonstration.} Appliquons à $f$ la prop. 5.3, en prenant pour $m$ le plus grand  entier tel que $a_m \neq 0$; on obtient $a_0 \geqslant2|a_m| \geqslant 2$ puisque $a_m$ est un entier non nul. Si $a_0=2$, la prop. 5.3 montre que  $f$ est divisible par, soit  $g_m^+= t^{-m}+2+t^m$, soit  $g_m^-=-t^{-m}+2-t^m$. Le terme non nul de plus haut degré du polynôme de Laurent $f/g_m^+$ (resp. $ f/g_m^-$) est nécessairement une constante, et cette constante est égale à 1
 puisque les coefficients de 1 dans $f, g_m^+, g_m^-$ sont les mêmes. On a donc, soit $f=g_m^+$, soit  $f = g_m^-$.

\medskip

\n {\bf 5.5. Applications aux $S$-caractères.}

\medskip
\n {\bf Proposition 5.6.}  {\it Le groupe $\U$ n'a aucun $S$-caractère $\neq1$.} 

\smallskip
\n {\it Démonstration.} Soit $f$ un $S$-caractère de U, distinct de 1. C'est un polynôme de Laurent $\sum a_nt^n$ positif, à coefficients entiers. Comme $\<f,1\> = 1$,  on a $a_0= 1$, ce qui contredit le cor. 5.4. 

\medskip

\n {\bf Proposition 5.7.} (= prop. 4.6) {\it Un tore n'a aucun $S$-caractère $\neq 1$.}

\smallskip
\n {\it Démonstration.} Soit $G$ un tore de dimension $n$. C'est un produit de $n$ copies de $\U$. Le cas $n=0$ est  trivial. Supposons $n \geqslant 1$. Soit $f$ un $S$-caractère de $G$. 
  Les caractères irréductibles de $G$ sont les éléments
du groupe  $ X = \Hom(G,\U)$. Notons $X$ additivement; c'est un groupe abélien libre de rang $n$; si $x\in X$, soit $\chi_x$ le caractère correspondant ; on a $\chi_{x+x'}= \chi_x\chi_{x'}$. Soit $n_x$ le coefficient de $\chi_x$ dans $f$. Soit $A$ l'ensemble des 
$x$ tels que $n_x\neq0$; c'est un ensemble fini, qui contient $0$ car $n_0=1$ puisque $f$ est un $S$-caractère. On a : $$ f = \sum_{a\in A} n_a\chi_a \quad {\rm et} \quad  n_0 = 1.$$ Soit $Y = \Hom(\U,G)$.  L'application
  $$ \Hom(\U,G) \times \Hom(G,\U) \to \Hom(\U,\U) = \Z,$$

\n permet d'identifier $Y$  au $\Z$-dual de $X$. Soit $y\in Y$
tel que $\<a,y\> \neq \<a',y\>$ pour tout couple d'éléments distincts $a,a'$ de
$A$; un tel élément existe car $A$ est fini. L'homomorphisme $y : \U \to G$ définit un caractère $f_y = f \circ y$
de $\U$. Si l'on identifie les caractères de $\U$ à des polynômes 
 de Laurent, on a
$$ f_y(t) =  \sum_{a\in A} n_at^{\<a,y\>}, \ {\rm avec} \ \  n_a \in \Z.$$
Vu le choix de $y$, on a $\<a,y\> = 0$ seulement si $a = 0$; le 
 terme constant de $f_y(t)$ est donc 1. Comme $f_y$ est positif, c'est un $S$-caractère de $\U$. D'après la prop. 5.6, on a $f_y=1$, ce qui signifie que $n_a$ = 0 pour tout  $a \neq 0$. D'où $f=1$.

\medskip

\n {\bf Proposition 5.8.} (= prop. 4.7) {\it Tout $S$-caractère du groupe $\SU(2,\C)$ est le carré d'un caractère irréductible.}

\smallskip 
\n {\it Démonstration}. Soit $G= \SU(2,\C)$, et soit $\chi$ un $S$-caractère de $G$. La restriction de $\chi$ au tore  maximal standard  U de $G$ est  un polynôme de Laurent
  $f = \sum_{n\in \Z} a_n t^n$. La moyenne  $\<\chi,1\>_G$ de  $\chi$ sur $G$ est donnée par un cas particulier élémentaire de la formule d'intégration de H. Weyl  ([Bo 82],  $\S$ 6.2, cor. 2 au th.1):
\medskip

\hspace{10mm} $\<\chi,1\>_G = \frac{1}{2} \<F,1\>_{\U}$, \quad {\rm où}  \ $F(t) = f(t)(-t^{-2}+2-t^2).$

\smallskip

Notons $b_i$ les coefficients  de $F(t)$; on a $\<\chi,1\>_G = b_0/2$; comme $\chi$ est un $S$-caractère, on a   $b_0=2$. D'autre part, $F(t)$ est produit de deux fonctions positives, donc est positif; il n'est pas constant. Le cor. 5.4, appliqué à $F$, montre qu' il existe un entier $m > 0$ tel que $F(t)$ soit égal à $t^{-m}+2+t^m$ ou à $-t^{-m}+2-t^m$. Comme $F(1)=F(-1)=0$, cela élimine $t^{-m}+2+t^m$, et montre que $m$ est pair. Si l'on pose $n=m/2$, on a donc
   $ F = -t^{-2n} + 2 - t^{2n}= - (t^n - t^{-n})^2.$
   Comme $-t^{-2}+2-t^2 = - (t-t^{-1})^2$, on en déduit
   
   \smallskip
   
  \hspace{15mm} $f(t) = g_n(t)^2$ \ où \ $g_n(t) =(t^n-t^{-n})/(t-t^{-1})$.
   
   \smallskip
 Soit $E_n$ une représentation irréductible de dimension $n$ de $G$. La restriction de $E_n$ au tore maximal U a pour caractère $g_n(t)$. Cela démontre la proposition.
   
   \bigskip

\n {\bf $\S$ 6. Nombre de zéros.}

\smallskip

  Revenons à la situation du $\S$ 1, où $G$ est un groupe algébrique
  linéaire sur un corps algébriquement clos $k$ de caractéristique $0$.
  
  Si $\chi$ est un caractère irréductible de $G$ de degré $>1$, nous avons vu
  qu'il existe au moins un zéro de $G$ d'ordre fini. On aimerait savoir s'il y en a seulement un nombre fini à conjugaison près et comment on peut les  
  déterminer. La question se pose, non seulement pour les caractères irréductibles, mais pour tous les caractères virtuels. Voici deux exemples:
  
  \smallskip
  (6.1) $G = \SL_2$, $\chi$ caractère virtuel $\neq 0$.
  
     \smallskip Si l'on écrit $\chi$ comme polynôme de Laurent, les zéros de $\chi$ sur le
    tore maximal standard sont ceux de ce polynôme, qui n'est pas $0$ par hypothèse. Ils sont donc en nombre fini.
    
      \smallskip(6.2) $G = \SL_2 \times \SL_2$, $\chi$ est le caractère de la représentation
    de degré 4  de $G$ qui est somme directe des représentations naturelles 
  de degré 2 de chaque facteur. 
  
  \smallskip  Un élément $t$ du tore maximal standard $T = \Gm \times \Gm$ de $G$ s'écrit de façon unique
  comme $t = (s_x, s_y)$, où $s_x = \begin{pmatrix} x & 0 \cr 0 & x^{-1} \end{pmatrix}$ et $s_y= \begin{pmatrix} y & 0 \cr 0 &y^{-1} \end{pmatrix}$.

  \smallskip  On a $\chi(t) = x + x^{-1}+y+y^{-1}$. Si $\epsilon$ est une racine de l'unité
  d'ordre $d$, et si $x=\epsilon, y=-\epsilon$, on a 
  $\chi(t)= 0$.  On obtient ainsi, pour tout $d$, un zéro de $\chi$ d'ordre $d$ si $d$ est pair, et d'ordre $2d$ si $d$ est impair. Ces zéros sont contenus dans $T_1 \!\cdot\! z$, où
  $T_1$ est la diagonale de $\Gm\times \Gm$ et  $z = (1,-1) \in T$. C'est là un cas particulier d'un résultat
  conjecturé par Lang et démontré par M. Laurent ([La 84]):
  
  \medskip
  
  \n {\bf Théorème 6.3}. {\it Soit $T$ un tore, et soit $f$ un élément de l'algèbre affine de $T$.
  Il existe un nombre fini de sous-tores $T_i$ de $T$, et des éléments
  $z_i$ de $T$ d'ordre fini, ayant les propriétés suivantes} :
  
  (i) {\it $f(x)=0$  pour tout $ x \in T_i\!\cdot\!z_i$}.
  
  (ii) {\it Tout élément $x$ de $T$ d'ordre fini tel que $f(x)=0$ est contenu
  dans l'un des $T_i\!\cdot\!z_i$.}
    
    \n [Noter que les éléments d'ordre fini de  $T_i\!\cdot\!z_i$ sont ceux qui s'écrivent $t_iz_i$ avec $t_i$ d'ordre fini dans $T_i$. Il y en a une infinité si $\dim T_i > 0$.]

  \smallskip Ce théorème s'applique à tout $f$ qui est un caractère virtuel de $T$.

\medskip
Pour utiliser le th.6.3 de façon effective, on doit déterminer les $(T_i,z_i)$, ce qui n'est pas facile lorsque $\dim T \geqslant 2$. Lorsque $\dim T  = 2$,  F. Beukers et C.J. Smyth ([BSm 00]) donnent une méthode explicite (voir $\S$ 7 ci-dessous). Dans le cas général, on a  des majorations du nombre de couples $(T_i,z_i)$, cf. [ASm 12]. 

   \medskip
\n {\bf $\S7$.  Exemple : éléments  $g$ de $\sf{G_2}$ d'ordre fini tels que la trace de $\Ad (g)$  soit $0$}

\smallskip
  Nous allons appliquer la méthode de [BSm 00] au caractère étudié au $\S$ 2.11: le groupe $G$ est un $k$-groupe simple de type $\sf{G}$$_2$, et $\chi$ est le caractère de sa représentation adjointe. Voici le résultat :

 \smallskip
 
 \n {\bf Théorème 7.1.} {\it Soit $n$ un entier $>0$. Pour qu'il existe $g\in G(k)$ d'ordre égal à $n$ tel que $\chi(g) = 0$, il faut et il suffit que $n = 7, \ 8, \ 15$ ou $42$.}
 
  \smallskip
 Soit $T$ un tore maximal de $G$ défini sur $k$. Comme tout élément semi-simple de $G$ est conjugué à un élément de $T$, et que les éléments d'ordre fini sont semi-simples, on voit qu'il suffit de démontrer le th. 7.1 pour les $g$ qui appartiennent à $T$.
  
 On identifie $T$ à $\Gm \times \Gm$ de telle sorte que les éléments $(1,0)$ et $(0,1)$ de $\Hom(T,\Gm)$  forment une base des racines de $G$. Si $t=(x,y)$ est un élément de $T$, on a alors :
 
  \smallskip
\ $\chi(t) = 2 + f(x,y) + f(x^{-1},y^{-1})$,
 
  \smallskip
 
 \n où \ $f(x,y) = x+ y + yx + yx^2+yx^3 + y^2x^3$.
 
 \smallskip
\n   {\small [Cela résulte de ce que les racines positives sont 
  $\alpha, \beta, \beta+\alpha, \beta+2\alpha, \beta + 3\alpha, 2\beta + 3\alpha$, 
  
  \n avec les notations habituelles ($\alpha$ racine courte et $\beta$ racine longue).]}
  
  \medskip
 
 \n Il est commode d'éliminer les dénominateurs, pour avoir des polynômes en $x,y$. Pour cela, on pose
 
 \smallskip
 $H(x,y) = y^2x^3 \chi(t) $,
 
 \smallskip
 \n et l'on a :
 
  \smallskip
 
$H(x,y)= y^4x^6+y^3(x^6+x^5+x^4+x^3) +y^2(x^4+2x^3+x^2)+y(x^3+x^2+x)+y+1$.

\medskip

\n Le théorème 7.1 est équivalent à :

\smallskip

\n {\bf Théorème 7.2}. (a) {\it Si $n = 7, 8, 15$ ou $42$, il existe des racines de l'unité $x$ et $y$ telles que $H(x,y)=0$
et que $t=(x,y) $ soit d'ordre $n$.}

(b) {\it Inversement, si $x$ et $y$ sont des racines de l'unité telles que $H(x,y) = 0$, alors l'élément $t = (x,y)$ est d'ordre $7, 8 , 15$ ou $42$.}

\medskip

\n {\it Démonstration de} (a). 

Pour tout $n >0$, notons $z_n$ une racine primitive $n$-ième de l'unité. Montrons  que, si $n = 7, 8,15$, on a $H(z_n,z_n^3) = 0$, et aussi $H(z_n, z_n^{11})=0$ pour $n=42$, ce qui prouvera (a). 

\smallskip

Pour cela, on considère le polynôme

 \smallskip
 
\n $H(x,x^3) = x^{18} +x^{15} + x^{14}+x^{13}+x^{12}+x^{10}+2x^9+x^8+x^6+x^5+x^4+x^3+1$,

 \smallskip
 \n et l'on constate qu'il se factorise de la façon suivante:

\smallskip
$H(x, x^3) = \Phi_7(x)\Phi_8(x)\Phi_{15}(x)$, 

\smallskip

\n où $\Phi_d(x)$ désigne le $d$-ième polynôme cyclotomique:

\smallskip
$\Phi_7(x) \ = x^6+x^5+x^4+x^3+x^2+x+1,$

 $\Phi_8(x) \ =x^4+1,$
 
$ \Phi_{15}(x)=x^8-x^7+x^5-x^4+x^3-x+1.$

\medskip

On a donc $H(x,x^3)=0$ si et seulement si  $x$  est une racine de l'unité d'ordre 7, 8 ou 15.

\smallskip

\n {\it Remarque.} C'est Deligne (lettre du 29/1/2024) qui m'a signalé les propriétés remarquables du sous-groupe de $T$
formé des éléments de la forme $(x,x^3)$, ainsi que la factorisation de $H(x,x^3)$.

\medskip

\n 
On a d'autre part

\smallskip
$H(x,x^{11})= x^{50}+x^{39}+x^{38}+x^{37}+x^{36}+x^{26}+2x^{25}+x^{24}+x^{14}+x^{13}+x^{12}+x^{11}+ 1$,

\n qui se factorise en

\smallskip
$H(x,x^{11}) = \Phi_8(x)\Phi_{42}(x)F(x)$,

\smallskip

\n où $\Phi_{42}(x) = x^{12}+x^{11}-x^9-x^8+x^6-x^4-x^3+x+1$, et :

\smallskip
\
$F(x)=x^{34}-x^{33}+x^{32}-x^{30}+x^{29}-x^{28}+x^{27}+x^{24}+x^{21}+x^{18}-x^{17}+x^{16}+x^{13}+x^{10}+x^7-x^6+x^5-x^4+x^2-x+1$.

\smallskip\n Cela complète la démonstration de la partie (a) du théorème.

\bigskip

\n {\it Démonstration de} (b).

\smallskip

Nous utilisons une méthode qui m'a été suggérée par C.J. Smyth, et qui  est exposée dans
[BSm 02], \S3.4. 

\smallskip

Définissons les sept polynômes en $x,y$ :

\smallskip
$H_1=H(x,-y), H_2= H(-x,y), H_3 =H(-x,-y), H_4=H(x^2,y^2), $

$H_5 = H(x^2,-y^2), H_6 =H(-x^2,y^2), H_7=H(-x^2,-y^2)$.

\smallskip
D'après {\it loc.cit.}, pour tout couple $t=(x,y)$ de racines de l'unité tel que $H(x,y)=0$, il existe un indice $i \in \{1, 2,..., 7\}$ et un seul tel que $H_i(x,y)=0.$ 

\smallskip Pour tout $i$, il nous faut donc déterminer quels sont les couples $(x,y)$ de  racines 
de l'unité tels que $H(x,y)= 0 = H_i(x,y)$.

  Cela peut se faire en éliminant la variable $x$ en formant le
  résultant de $H, H_i$, considérés comme polynômes en $x$, à coefficients dans le  corps $\Q(y)$. On obtient ainsi
  un polynôme $S_i(y)$. On factorise $S_i(y)$ en produit de polynômes irréductibles, et l'on ne garde que les facteurs qui sont des polynômes cyclotomiques: sur [PARI], cela se fait au moyen de la  commande {\it polcyclofactors}.\footnote{C'est aussi cette commande qui m'a permis de détecter la divisibilité de $H(x,x^{11})$ par $\Phi_{42}$.}
  
  Les racines de ces facteurs constituent la liste des $y$ possibles. En permutant les rôles de $x,y$ on obtient de même un polynôme $R_i(x)$, d'où une liste des  $x$ possibles. Si $y_0,x_0$ appartiennent à ces listes, on calcule $H(x_0,y_0)$
et l'on voit s'il est ou non zéro (comme $H(x_0,y_0)$ est un entier d'un corps cyclotomique dont tous les conjugués
sont de module au plus 14, un calcul approché est suffisant pour déceler zéro).

\medskip
Le tableau suivant donne les résultats de ces calculs\footnote{Faits avec l'aide de C.J. Smyth, et de [PARI].}, pour chaque valeur de $i = 1,\dots, 7$, avec les conventions suivantes:

 $z_d$, pour $d = 2,3,5,7,8,15,42$ : une racine primitive $d$-ième de l'unité;

 $ R_i(x)^{cycl} $  : l'ensemble des polynômes cyclotomiques qui divisent $R_i(x)$;
 
   $S_i(y)^{cycl}$  :  l'ensemble des polynômes cyclotomiques qui divisent $S_i(y)$;
   
 couples : les couples possibles $(x,y)$ pour chaque valeur de $i$; 
   
  ordre du couple: les ordres des éléments $t=(x,y)$ ainsi obtenus, i.e. le ppcm des ordres de $x$ et $y$.

       $$\begin{array}{|c|c|c|c|c|l}
\cline{1-5}
\hbox{$i$} &  R_i(x)^{cycl} &  S_i(y)^{cycl}  &  {\rm couples}  & {\rm ordre \ de} \ t= (x,y) &\\
\cline{1-5}
1& \Phi_2  \ \ \Phi_4 & \Phi_8 &(z_2,z_8) \ \ (z_8^6,z_8) &8 \ \ \ 8&\\
\cline{1-5}
2& \Phi_8 & \Phi_2 \ \ \Phi_4 & (z_8,z_2) \ \ (z_8,z_8^2) & 8 \ \ \ 8 &\\
\cline{1-5}
3& \Phi_8 & \Phi_8 & (z_8,z_8) \ \ (z_8,z_8^3) & 8 \  \ \ 8&\\
\cline{1-5}
4& \Phi_3 \ \ \Phi_7 \  \ \Phi_{15}& \Phi_5 \ \ \Phi_7 & (z_3,z_5) \ \ (z_7,z_7) \ \ (z_7,z_7^3) \ \ (z_{15}, z_{15}^3) \ \ (z_{15},z_{15}^9) & \ 15 \ \ \ 7 \ \ \ 7 \ \ \ 15 \ \  \ 15 &\\
\cline{1-5}
5 & \Phi_7 & \Phi_2 \ \ \Phi_{42} & (z_{42}^{18},z_{42}) & 42 &\\
\cline{1-5}
6 & \Phi_{42} & \Phi_3 & (z_{42},z_{42}^{28}) & 42 &\\ 
\cline{1-5} 
7 & \Phi_{42} & \Phi_2 \ \ \Phi_{42} & (z_{42},z_{42}^{11}) & 42 &\\
\cline{1-5}

\end{array}$$

\medskip

On constate que les seuls ordres possibles sont ceux du théorème 7.2.
 
\newpage

\hspace{45mm}{\bf Références}

\vspace{1cm}

\n [ASm 12] I. Aliev \& C.J. Smyth, {\it Solving algebraic equations in roots of unity}, Forum Math. {\bf24} (2012), 641--665.

 \n [BKZ 19] Y. Berkovich, L.S. Kazarin \& E. Zhmud, {\it Characters of Finite Groups},  vol. 2, de Gruyter, Berlin, 2019.

\n [Bo 72] N. Bourbaki, {\it Groupes et algèbres de Lie}, chap. 2--3, Hermann, Paris, 1972.

\n [Bo 68] \ ---------, {\it Groupes et algèbres de Lie}, chap. 4--6, Hermann, Paris, 1968.

\n [Bo 75] \ ---------, {\it Groupes et algèbres de Lie}, chap. 7--8, Hermann, Paris, 1975.

\n [Bo 82] \ ---------, {\it Groupes et algèbres de Lie}, chap. 9, Masson, Paris, 1982.

\n [Bo 23] \ ---------, {\it Théories spectrales}, chap. 3--5, Springer, 2023.

\n [BSe 64] A. Borel \& J-P. Serre, {\it Théorèmes de finitude en cohomologie galoisienne}, Comm. Math. Helv. {\bf39} (1964), 111--164 (= A. Borel, Oe. II, n°64).

\n [BSm 00] F. Beukers \& C.J. Smyth, {\it Cyclotomic points on curves}, Number Theory for the Millennium I (Urbana, IL, 2000), 67--85, A. K. Peters, Natick, MA, 2002.

\n [BJM 24] T. Breuer, M. Joswig \& G. Malle, {\it Zeros of S-characters},
arXiv 2408.16785.

\n [Bu 03] W. Burnside, {\it On an arithmetical theorem connected with roots of unity and its application to group characteristics}, Proc. L.M.S. {\bf1} (1903), 112--116.

\n [CC 92] P.J. Cameron \& A.M. Cohen, {\it On the number of fixed point free elements in a permutation group}, Discrete Math., {\bf106/107} (1992), 135--138.

\n [Ch 46] C. Chevalley, {\it Theory of Lie Groups}, Princeton U. Press, 1946.
  
  \n [Dy 52] E.B. Dynkin, {\it Sous-algèbres semi-simples des algèbres de Lie semi-simples} (en russe), Mat. Sbornik {\bf30} (1952), 349--462 (=  Selected Papers, 175--308).
     
\n [EKV 09] A.G. Elashvili, V.G. Kac \& E.B. Vinberg, {\it On exceptional nilpotents in semisimple Lie algebras}, J. Lie Theory {\bf19} (2009), 371--390.

\n[FKS 81] B. Fein, W.M. Kantor \& M. Schacher, {\it Relative Brauer Group} II, J. Crelle, {\bf328} (1981), 39--57.

\n [Ja 03] J.C. Jantzen, {\it Representations of algebraic groups}, second edition, A.M.S. SURV {\bf107}, 2003.

\n [Jo 72] C. Jordan, {\it Recherches sur les substitutions}, J. Liouville {\bf17} (1872), 351--367 (= Oe. I, n°52).

\n [Ka 81] V.G. Kac, {\it Simple Lie groups and the Legendre symbol}, LNM {\bf848} (1981), 110--123.

\n [Ko 59] B. Kostant, {\it The principal $3$-dimensional subgroup and the Betti numbers of a complex simple group}, Amer. J. Math. {\bf81} (1959), 90--98 (= Coll. Papers, vol. I, n° 11).

\n [La 84] M. Laurent, {\it \'{E}quations diophantiennes exponentielles}, Inv. math. {\bf78} (1984), 299--327.

\n [Ma 24] G. Malle,  {\it Zeros of characters}, arXiv 2408.16785.

\n [MNO 00] G. Malle, G. Navarro \& J.B. Olsson, {\it Zeros of characters of finite groups}, J. Group Theory {\bf3} (2000), 353--366.

\n [PARI] PARIS/GP  http://pari.math.u-bordeaux.fr/.

\n [Se 93]  J-P. Serre {\it Gèbres}, Enseign. Math. {\bf39} (1993), 33--85 (= Oe. IV, n°160).

\n [Se 98] \ ---------, {\it Représentations linéaires des groupes finis}, deuxième édition refondue, Hermann, Paris, 1998.
 
  \n [Se 04]  \ ---------, {\it  On the values of characters of compact Lie groups}, Oberwolfach Reports {\bf1} (2004), 666--667.
         
     \n [Si 50] J. de Siebenthal, {\it Sur certains sous-groupes de rang $1$ des groupes de Lie clos}, C.R.A.S. {\bf230} (1950), 910--912.

\n [St 68] R. Steinberg, {\it Lectures on Chevalley Groups}, 1968 Yale Notes, prepared by J. Faulkner and R. Wilson, ULECT {\bf66}, AMS, 2016.

\n [We 26] H. Weyl, {\it Theorie der Darstellung kontinuierlicher halb-einfacher Gruppen durch lineare Transformationen.} III, Math. Zeitschrift, {\bf24} (1926), 377--395.

\n [Zh 95] E.M. Zhmud,  {\it On one type of nonnegative generalized characters of a finite group}, Ukrainian Math. J., {\bf47} (1995), 1526--1540.

\bigskip

\n Collège de France, 11 Place Marcelin Berthelot, 75231 Paris

\end{document}